# Evidencing missing resources of the documentational approach to didactics. Towards ten programs of research / development for enriching this approach[1]

Luc Trouche, French Institute of Education, ENS de Lyon, France

**Abstract:** This chapter proposes a view from inside the DAD, starting from determining some essential resources missing of DAD, to proposing 10 programs of research/development for developing it. It could be considered as a follow-up of Chapter 1, where Ghislaine Gueudet situates the current state of DAD in looking back to its origin: chapter 12 proposes a possible future of this approach in analyzing its current state. It determines the missing resources of DAD in questioning current and past PhD students who have anchored their research in DAD. What did/do they learn in using DAD as a main theoretical resource; to which extent did/do they estimate that they have enriched DAD by their own work? Which are, according to them, the still missing resources of DAD? Which of these resources should be developed by DAD from itself, and/or in co-working with other theoretical framework? From this inquiry, this chapter proposes ten perspectives of research, aiming to develop some theoretical blind points of DAD, or to develop some methodological tools, or to deepen the cultural/social aspects of DAD in questioning the naming systems used by teachers when interacting with resources. This chapter echoes actually different perspectives of research already present, as promising germs, in previous chapters of the book.

**Keywords:** Documentational approach to didactics, missing resources, research program, resource system, theoretical networking,

**Introduction**

This chapter is dedicated to the future of the Documentational Approach to Didactics (DAD). Looking towards the future usually starts with 'looking back', and questions such as the following arise: where do we come from, which was the path we followed for arriving here, which issues did we have to face, and which were the 'missing' resources in our view that we need to advance DAD?

My first book (Trouche et al. 1998) was written with my 37 students, while I was still teaching mathematics in a secondary school: *Experiencing and proving. Experimenting mathematics in schools with symbolic calculators, 38 variations on a given topic*. It was a great experience, for me, to learn from my students regarding my own teaching, the potential of the mathematical problems we had faced together, and the 'missing' resources in the mathematical environment of the classroom. Twenty years after, I will situate my talk in the same spirit, drawing from interactions with students and young researchers having used DAD, for retrospective and prospective reflections.

---

[1] This chapter originates from a lecture given to the Re(s)sources 2018 International Conference. Video in English, with French subtitles, available at http://video.ens-lyon.fr/ife/2018/2018-05-30_009_Ressources2018_Luc_Trouche_v1.mp4

These reflections result from ten programs of research, the number 10 resonating with the 10 years of DAD. Actually, I wished for a reasonable number, and 10 appeared as a good compromise between a too small number of huge programs, and a too big number of scattered programs. Nevertheless, I was thinking of these programs in a kind of 'free attitude' mood: I wrote this chapter shortly before retiring as professor emeritus, having in mind that I will not be in charge of coordinating these programs (but hoping to participate in some of them!). Also, I was thinking of these programs more in terms of necessity, than in terms of feasibility, in an essential perspective of transmission to the community that I have tried to contribute to during the previous 10 years. When, in the following, I use the "we", it will designate this community, i.e. people studying teachers' work with resources under, at least partially, the umbrella of DAD.

After this introduction, this chapter is structured in 7 sections: after this introductory section, in section 1, I present the way of investigation I have used for detecting' DAD existing and missing resources, and inferring the corresponding research programs; section 2 evidences some main DAD resources; section 3 focuses on some missing *theoretical* resources; section 4 on some missing *methodological* ones; section 5 questions the necessity of an extension/expansion of the theory; section 6 turns back to what is grounding each theory in terms of: history, culture, and, finally, words; and section 7 discusses the processes to be continued, carried on, or engaged.

**1. Learning from those who have appropriated and developed DAD**

In this section, I briefly situate the roots of my experience as a researcher, then I present the way I have gathered the data for designing this chapter.

*1.1 Rooting a personal experience*

From my own experiences, starting when I was teaching (and learning), and designed or used resources for this purpose. I have already cited the work done with my secondary school students. As a researcher, I have already tried to describe my own intellectual trajectory (Trouche 2009), made of encounters with: teachers, teacher educators and researchers; projects of research; contrasted cultural situations; and students.

Tracing the encounters with researchers is often easy, as these encounters usually produce scientific resources (papers or books): developing an instrumental approach to didactics, with Dominique Guin (Guin and Trouche 1998), then with Kenneth Ruthven (Guin et al. 2005); developing the notion of instrumental orchestration with Paul Drijvers (Drijvers and Trouche 2008), developing the documentational approach to didactics with Ghislaine Gueudet (Gueudet and Trouche 2009), associating then Birgit Pepin (Gueudet et al. 2012); developing a wider reflection on tools and mathematics with John Monaghan (Monaghan et al. 2016). Besides, some of the publications mentioned before were collective books; acting as an editor provided me the opportunity to meet the authors participating to the present book.

Tracing the encounters with the diversity of actors of mathematics education, teachers, teacher educators, researchers, students, engineers, is, generally, not so easy, because these encounters occurred in a variety of organizations and contexts. I keep mainly in mind my experience as the director of the Institute of Research on Mathematics Teaching of Montpellier (Trouche 2005a), as the president of the French Commission on Mathematics Teaching (Trouche 2017) and as a member of the EducTice Team of the French Institute of Education. Three research projects had also a major role for developing DAD, and crossing it with other approaches: a regional project, PREMaTT[2]; a national project, ReVEA[3]; and a European project, MC2[4].

I also learnt a lot from my experiences in scientific stays abroad: In Brazil (UFPE, Recife), where I deepened the idea of *webdocuments* (Bellemain and Trouche 2016), in China (ECNU, Shanghai), where I better understood the importance of teachers' collective work (Wang et al. 2018), in Mexico (Cinvestav, Mexico DF), where I better understood the importance of socio-cultural approaches (Radford 2008), and in Senegal, where I discovered the missing teaching resources of developing countries (Sokhna and Trouche 2016). A symposium, co-organized with Janine Remillard and Hendrik Van Steenbrugge in the frame of the second International Conference on Mathematics Textbooks (2017), *Teacher-Resource Use around the World*, offered a good opportunity for crossing these experiences.

Finally, I believe that my main learning sources were the PhD students themselves, appropriating, using and developing DAD. We really know indeed what we try to appropriate, what we use, what we experience, what we develop in a creative way. What is true for any artefact is all the more true for a theoretical framework. Till now, 10 students have defended their PhD, considering DAD as part of their theoretical framework[5]. About 10 PhDs are in progress in the same frame. And a number of PhDs and post-doc students are more or less exposed to DAD through scientific stays or research projects. Obviously, the different elements rooting my own experience are not isolated: PhD, research projects and international collaboration are linked, for example via PhDs co-supervision. But my choice, for conceiving this chapter, was to privilege the direct feedback of students, keeping in mind the fundamental dialectic between learning and teaching. I like, in this perspective, the Chinese translation of "Teaching", which is 教学 (Jiàoxué: a concatenation of two characters meaning respectively "teaching" and "learning").

*1.2 Method for gathering data grounding this chapter*

I collected students' feedback by two means: a seminar and a questionnaire.

The monthly 'Resources seminar' was organized in the French Institute of Education, from September 2017 to June 2018, gathering about 15 PhD and post-doc students[6]. Each session of this seminar gave the opportunity to two students for presenting and connecting their work. It allowed for questioning the main concepts of DAD and for surfacing evidence of some missing resources.

The questionnaire was launched in December 2017, the answers being expected for March 2018. It was sent to students or young researchers using, or having used, DAD in their PhDs. I did not aim to reach all the people filling these conditions, but only the students and young researchers that I knew closely enough, with respect to their general research work. I got 29 answers from the 32 people solicited (see the corresponding names in the final acknowledgements, § 8). In the letter accompanying the questionnaire, I motivated it by the necessity, for me, to prepare a lecture for the Re(s)sources International Conference, and by my intention to base my lecture on their reflective view of DAD, for benefiting to the development of the research community itself. Then I asked a limited number of questions, and left it up to people to more or less develop their answers.

---

[2] PREMaTT: thinking the resources of mathematics teachers in a time of transitions (http://ife.ens-lyon.fr/ife/recherche/groupes-de-travail/prematt)

[3] ReVEA: Living resources for learning and teaching (https://www.anr-revea.fr)

[4] MC2: Mathematical creativity squared (http://www.mc2-project.eu/)

[5] See the list of the corresponding PhD here: https://ens-lyon.academia.edu/enslyonacademiaedu/PhD

[6] http://eductice.ens-lyon.fr/EducTice/seminaires/ressource-2017-2018

The questions are:

1. You have used instrumental approach (IAD) and/or documentational approach to didactics (DAD). To what extent did these frames support your research?
2. In using DAD, did you feel that there were resources missing? In which circumstances? What new resources (theoretical as well as methodological) should this approach develop?
3. For designing missing resources, which possible theoretical networking? Which new research programs to be launched?

The first question mentions the instrumental approach to didactics, taking into account that: some students finished their PhD before the introduction of DAD (e.g. Sokhna 2006); for them, DAD appeared as a natural follow-up of their previous research. Some other students (e.g. Lucena et al. 2016) are using the notion of instrumental orchestration, situating themselves somewhere in between IAD and DAD. Actually, IAD acted as an incubator of DAD (see Gueudet, Chapter 1).

I started my writing from the answers given to these questions, and from the interactions having occurred in the monthly Resources seminar. I have tried of course to take into account the whole set of answers. Then I crossed, gathered and structured them according to my own experience (§ 1.1 above…). Doing this, I was aware of some limitations coming from the methodology itself: some students may have thought that their difficulties come from their own lack, and not from the theoretical framework itself. This is the reason why the formal and informal interactions with these students were important, for complementing some points in the students' answers. This work gave matter to ten programs of research and/or development.

## 2. The productive and constructive aspects of a theoretical framework

In this section, I present essentially the answers to the first question above, distinguishing three parts: (1) the way DAD supports the research; (2) the resources already developed for/by this approach; (3) the main missing resource.

*2.1 To what extent does DAD support your research?*

The contributions of the DAD could be classified in three main categories, *understanding*, *designing* and *rethinking*.

*Understanding* concerns mainly the resources in their *diversity* and their *complexity*, and for their *potential* for resourcing teachers work; but also as critical *interfaces* between teachers and students, or between *individual* and *collective* work, or between designers and users and as markers; or *witnesses* of teachers' professional development.

Understanding concerns then:

- The *potential* of various things for resourcing teachers' work;
- The *complexity* of the integration processes;
- The *interactions* between teachers and *students*;
- The interactions between *individual* and *collective* teachers' work;
- The gap between designers' *expectations* and teachers' *uses*;
- The metamorphosis from *prospective teachers* to *teachers at work*;
- Teachers' professional development over a *long period*.

*Designing* concerns mainly a didactical engineering dimension of the research work, focusing on the role of digital resources for renewing tasks and supports to teachers' work.

Designing concerns then:

- Tasks integrating *digital resources*;
- *Digital supports* for helping teachers to design their own resources.

*Rethinking* concerns teacher education, and the nature itself of « what is teaching, what is a teacher », with a double component of designers and collective-reflective practitioners.

Rethinking concerns then:

- *Preservice as well as in-service teacher education*;
- Teachers themselves as resource *designers* and *collective-reflective practitioners*.

What I retain from these answers is the central role of the notion of 'resources' in the emergence of DAD, and more specifically, the role of digital resources; and the double potential of the approach for its practical component (supporting design) and its conceptual component (rethinking teaching and teacher education). This could be compared to the emergence of the instrumental approach to didactics (e.g. Guin and Trouche 1998), but the issue was, at this time, mainly the integration of one new artefact (in this case symbolic calculators). DAD proposes a holistic approach to teachers' work, taking into account the new universe of resources offered to teacher use, design and re-design.

*2.2 The new resources developed by/for DAD over the past 10 years*

What appears clearly during these last 10 years, and that is always the case when a new approach is proposed, is the flourishing process of creating concepts and names: first of all about the resources themselves, but also about the events occurring over the time of designing and using, about the knowledge guiding the process of design and use, about the collective sheltering these processes, and about the methodology of research, quite complex, to be constructed due to the diversity of times and places where the issues of teaching resources are addressed.

The answers to the questionnaire mention this creative process: being supported by a given frame goes on with contributing to the development of this frame, in the dynamics of each PhD:

- About resources: mother and daughter resources, complemented by the notions of structuring mother resources and oriented daughter resources; pivotal resources; meta-resources; proper, recycled and intermediary resources; block of resources; cycle of life of a given resource; resource systems for learning (resp. for teaching)

- About events: documentational incidents; didactical incident;
- About knowledge: didactical affinity, documentational identity; documentational expertise; information competencies, documentational experience and documentational trajectory;
- About collectives: mother, daughter collective, and hub collective (ReVEA project), documentation-working mate;
- About methodological design: methodological contract, reflective and inferred mapping of a resource system)

Of course, other conceptual developments occurred during the past 10 years (see Gueudet, Chapter 1), I just mention here the contributions evoked over the questionnaire answers. The awareness of this conceptual diversity leads naturally to the need for a map allowing us to circulate into this new field.

*2.3 Research program n°1: designing a DAD living multi-language glossary*

The need for a tool allowing one to master the conceptual field of DAD appears in a number of answers. Actually, a seed of such a tool had been developed at the emergence of the approach, by Ghislaine Gueudet, Rim Hammoud and Hussein Sabra (the two first PhD students situating their work in DAD) and me, on a website, mainly in French, dedicated to DAD (http://educmath.ens-lyon.fr/Educmath/recherche/approche_documentaire) and offering a glossary presenting a definition of 21 critical terms grounding DAD, supported by the first papers written in this period (see Figure 1).

> *Curriculum support:* Resource elaborated for a teaching purpose (see the notion of *curriculum material*, Remillard 2005).
>
> Remillard, J.T. (2005). Examining key concepts in research on teachers' use of mathematics curricula. *Review of Educational Research 75*(2), 211-246.

**Figure 1.** An example (our translation) from the glossary developed at the emergence of DAD

But this seed was not developed further, as it is always difficult to combine the development of a field, and the development of a map allowing us to explore this expending field …

A new initiative was more recently launched, in the frame of the ReVEA project: The AnA.doc platform (Alturkmani et al. to be published) gave access to a set of *situations* of teacher documentation work (see also § 4.4), and to a glossary integrating definitions of the main concepts used in the analyses (see Figure 2). Then the glossary evolves with the analyses of the situations themselves: new concepts appear, and the meaning of each of them could also evolve, according to the evolution of the theory itself. The problem, here, is that the AnA.doc platform was developed in the frame of a given project, in French, and it works like a prototype, for a rather small community[7].

---

[7] AnA.doc is available at https://www.anr-revea.fr/anadoc. A password is needed for accessing it, to be asked to the author of this chapter.

| **Situations side of AnA.doc** |
|---|
| Situation: Sophie (a middle school mathematics teacher), evokes her professional experience with resources.<br><br>Description: In a one-hour interview, in her school, we asked Sophie about the events over time that influenced her interactions with her **resources** for teaching. Due to her different roles (member of a teacher association, teacher educator…), the events linked to her **collective work** are particularly underlined.<br><br>From the information she gave, we tried to infer her **resource system** and her **documentational trajectory**. |
| **Glossary side of AnA.doc** |
| *Documentational trajectory*: this concept is currently developed in the frame of the Rocha's PhD. For now, a documentational trajectory is constituted by all the events influencing over time, the teacher documentation work, i.e. both her resource system and the associated knowledge (Rocha, 2016) |

**Figure 2.** An extract of the AnA.doc platform (our translation): each expression in bold (in the situations side) has a corresponding definition in the glossary side

Drawing from my experience of teaching in different contexts, I think that we need to combine the design of a glossary and the reflection of its instantiation in different languages. I realized the interest of such work in teaching, in English, in China (Figure 3): explaining a definition in a very different language provides opportunities to deepen the corresponding concept, to give examples and counter examples. The process of *denominating concepts* is essential for a development of each scientific field (Rousseau and Morvan 2000).

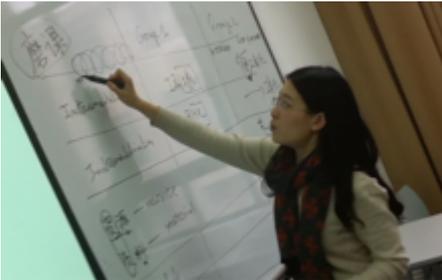

A twofold discussion:
- What could be a Chinese explanation, and, if possible, translation, of critical terms of DAD such as *instrumentation* or *instrumentalisation*?
- What could be an English explanation, and, if possible translation, of critical terms used for describing teachers' documentation work, such as 磨课 (Mó kè)?

This explanation needed (see picture right side) actually a lot of words and gestures[8]?

**Figure 3.** An episode of a master class, in ECNU (Shanghai): discussing the concepts and their translation

---

[8] MOKE is a method metaphor of refining a lesson as sharpening (a knife). It is a circulative-rise mode based on an iterative process: the lesson is designed and implemented a first time in the context of an open class (other teachers attending and observing this lesson), then discussed among the audience. After a refinement based on the discussion, a new version is implemented again in another class, etc.

Such work has existed in other emergent fields, for example in the field of TEL (Technology Enhanced Learning), with two meta-projects, intending to create an intellectual platform to support the conceptual and theoretical integration in this research area: a TEL Thesaurus and a TEL Dictionary. "Both tools are fully interdisciplinary, multilingual and takes into account the multicultural and epistemological roots of research on learning" (https://www.tel-thesaurus.net). Far from hiding the difficulty arising from translation processes, this platform profits from these issues for enriching the concepts at stake. For each language taken into account by the platform, a list of terms is proposed. An editor is in charge of each of them, and proposes for it: a definition, comments on its history, related terms, the translation issues, the disciplinary issues, and some key references. The editor is also in charge of accepting and managing the contributions proposing complementary or alternative views on each aspect of the "identity card" of this term in a given language.

**It leads to my first potential program of research/development: developing a DAD living and multi-language glossary.**

For me, it is really a condition for developing a scientific community, its concepts and methods. It is also a condition for each researcher to express his/her analyses and write papers in his/her own language, necessary both for the dissemination of a new theoretical frame, and its enrichment in encountering other cultures of teaching and of research.

**3. Back to basics, deepening the model of DAD itself, through four research programs**

I address in this section some theoretical issues related to DAD. Starting by the missing resources pointed out by the answers to the questionnaire, I propose afterwards four potential programs of research aiming to take into account these theoretical needs.

*3.1 The missing theoretical resources*

Ten years is a quite short period for the genesis of a theoretical frame, crossing different scientific fields (see Gueudet Chapter 1). In this initial phase, some fuzzy aspects could appear. In his answer, one PhD student pointed out some wavering, for example about the notion of document, often presented with an equation (different languages used, our translation in English):

- In Gueudet and Trouche (2010): a document = recombined resources + a scheme of utilization;
- In Trouche (2016): a document = combined resources + schemes of utilization;
- In Bellemain and Trouche (2016): a document = resources + a scheme of instrumented action;
- In Pepin et al. (2017): a document = the joint resources + their usages + knowledge guiding their usages

Is there one or several schemes developing with a given document? What means "utilization" in the expression "scheme of utilization" (and what is the balance between using and designing)? Are usages only *guided* by teacher's knowledge (and what about usages *producing* new knowledge?). This conceptual wavering reveals the underlying complexity of the processes at stake. I have classified answers to the questionnaire into four categories, that could be related to different sensitive points of the theoretical model:

- A category for *resources and resource systems*: what could be considered as a resource (gestures, languages, artefacts…)? What are the resource systems components? What about students' resource systems, and their interactions with teachers' resource systems?
- A category for the dialectic *schemes / situations*: how could they be simultaneously analyzed in their joint evolution?
- A category for the dialectic *between individual vs. collective documentation work of teachers*, their context and effects;
- A category for *the lenght of teachers' professional development* (dialectic between short vs. long geneses).

Drawing from these categories, I propose the following four potential research programs.

*3.2 Research program n°2: Modeling the structure and the evolution of teachers' resource systems*

This program takes into account the first category of students' propositions (§ 3.1), answering the questionnaire, about *resources* and *resource systems* (issues already addressed in Chapters 1, 2 and 8).

*About the notion of resource*, from the beginning of DAD, we rely on the Adler's productive metaphor: "I also argued for the verbalization of resource as 're-source' » (2012, p. 4). Doing this, we cultivated a kind of ambiguity, using the term resource sometimes for something having the *potential* to re-source teachers' activity, and sometimes for something already integrated in teachers' activity. The following quotation (Trouche et al. 2018b Online First) is interesting from this point of view:

> "Retaining [Adler's] point of view, DAD took into consideration a wide spectrum of 'resources' that have the potential to resource teacher activity (e.g., textbooks, digital resources, emails exchanged with colleagues, students' sheets), resources speaking to the teacher (Remillard 2005) and supporting her/his engagement".

This quotation is interesting, because it begins with the notion of "resource as a potential", then it goes on with the notion of a "resource interacting with a teacher", and finally ends by the notion of a "resource supporting a teacher".

Actually, when we use the expression "curriculum resources", we use it in a first meaning of "potential resources"; when we use the expression "a teacher's resource", we use it in the second meaning of a resource already *adopted* (which often means *adapted*, as a result of the instrumentalization processes) by him/her. Finally, we have to make clear that this second meaning is the one that DAD is retaining. This choice is coherent with the whole sentence of Adler:

> "I also argued for the verbalization of resource as 're-source'. In line with 'take-up', I posited that this discursive move shifts attention off resources per se and refocuses it on teachers working with resources, on teachers re-sourcing their practice (2012, p. 4).

It was also what we said in the seminal paper introducing this approach (my translation): "What the activity of a teacher encompasses, it is a set of resources" (Gueudet and Trouche, 2008, p. 7). With this idea, comes a direct consequence: "A resource is never isolated; it belongs to a set of resources" (Gueudet and Trouche 2009, p. 205). This essential distinction between resource per se and resources integrated into the activity of a given teacher opens a wide space for issues:

- If we consider the resources *per se*, how could we analyze their potential (epistemic, didactical), their variability, and their quality?
- If we consider the resources *of a given teacher*, how could we capture, and analyze their diversity, from gestures to artefacts (Salinas and Trouche 2018) and their role all along his/her documentation work?
- How could we capture and analyze the process from potential resources to resources integrated into teachers' practices? Hammoud (2012) has introduced the notions of *mother resources* and *daughter resources*, Alturkmani et al. (2018) the notions of *structuring mother resources*, and *oriented daughter resources*: some first steps for analyzing the dialectic between *potential and actual resources*.

*About the notion of resource system*. Considering a teacher' resource leads then to consider the set of resources it belongs to, that is *their resource system*. The choice of the expression resource *system* constitutes a strong statement, meaning that the set of resources that a teacher had appropriated is not messy, but is organized according to a given structure. Besides, from the beginning of DAD, a teacher's resource system is defined in an indirect manner: "The resource system of the teacher constitutes the 'resource' part of her documentation system (i.e. without the scheme part of the documents)" (Gueudet and Trouche, 2012a, p. 27). This definition leads to what I call a top-down perspective, starting from the analysis of a part of the documentation system (i.e. a document developed for a give purpose), for inferring the resources involved in this genesis. I plead for a more balanced point of view, combining this top-down perspective, and a bottom-up perspective (from the resource system to the documentation system). Considering the set of a teacher's resources as a *system* opens also a wide space of questions:

- From which point of view could we analyze this structure (from a location point of view, from the mathematics knowledge point of view, from a curricular point of view, from a didactical point of view…), and how can we articulate these different points of view?
- How could we characterize such a system? As an open dynamic system, or…? And which consequence should have such a characterization?
- Analyzing a given system leads to analyze specific points, and specific links. Some works have already (see Gueudet, chapter 2) pointed out typical resources, as pivotal, or meta, resources; we need certainly to deepen this classification;
- Could this analyze lead to evidence, for a given teacher, missing resources (and from which point of view?), as Chevallard and Cirade (2012) did, considering not the resources of a given teacher, but the resources of the profession itself?
- For a given teacher, we could distinguish different resource systems intersecting his/her own resource system: the classroom resource system (see Ruthven, chapter 1), his/her students' resource systems, and the resource systems of various collectives s/he participates to. Analyzing the *interfaces* (Trouche et al. 2018a Online first) between these resource systems could shed some light on the own teacher' resource system.

There are many ways for addressing the issues related to resource systems (analyzing their links to documentation systems, for example), but the joint study of the resources themselves and the system structuring them seems to be promising, and correspond to the first set of answers to the questionnaire. Therefore my second potential research program will be: **Modeling the structure and the evolution of teachers' resource systems**

*3.3 Research program n°3: Deepening the dialectics of schemes / situations of documentation work*

This program takes into account the second category of students' propositions (§ 3.1), answering the questionnaire, about the dialectic schemes / classes of situations. These propositions question the model grounding DAD (Figure 4), asking for further developments at three levels: the notion of scheme, the notion of situation, and the structure of documentation work (issues already addressed in Chapter 10).

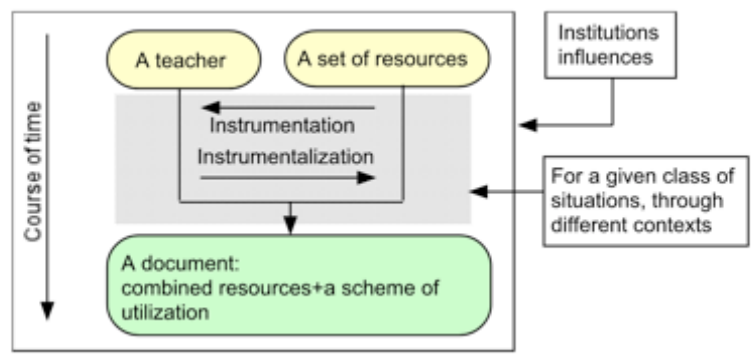

**Figure 4.** A representation of a documentational genesis (see Gueudet, Chapter 1)

The first level of reflection is related to the notion of scheme, which is certainly complex (Vergnaud 2009). But the notion of *utilization scheme* should be particularly questioned. If we go back to the instrumental approach to didactics, I had distinguished (Trouche 2005a, p. 150) *usage schemes and instrumented action scheme*, following the Rabardel's distinction:

- *Usage sche*mes related to 'secondary tasks'. These can be located at the level of elementary schemes (meaning they cannot be broken down into smaller units liable to meet an identifiable sub-goal), but it is by no means necessary: they can themselves be constituted as wholes articulating a set of elementary schemes. Their distinctive feature is that they are orientated towards secondary tasks corresponding to the specific actions and activities directly related to the artifact;

- *Instrument-mediated action schemes* which consist of wholes deriving their meaning from the global action which aims at operating transformations on the object of activity. These schemes incorporate usage schemes as constituents. Their distinctive feature is their relation to 'primary tasks'. They make up what Vygotsky called 'instrumental acts', which, due to the introduction of the instrument, involve a restructuring of the activity directed towards the subject's main goal […]. Usage schemes constitute specialized modules, which, in coordination with one another and also with other schemes, assimilate and mutually adapt in order to constitute instrument-mediated action scheme (Rabardel 2002, p. 83).

A scheme is developed in order to accomplish a given task, then associating the development of a document to an scheme of utilization is certainly reducer. Some works (for example Messaoui 2018) have begun to address this issue. I think that we could distinguish *usage scheme* (for example: storing a new resource in his/her own resource system) and *document action scheme*, oriented 'towards the subject's main goal' (for example: preparing a given lesson). It opens new questions, in the context of documentation work, towards a typology of usage schemes, and the possible decomposition of documentation action scheme in different usage schemes.

The second level of reflection is related to the DAD model (Figure 4) itself. It was very productive, allowing a lot of deep analyses of teachers' documentation work over the last 10 years. Pointing out the dialectical processes of instrumentation and instrumentalization allowed, on one side, to study the effects of resources on teacher's activity, and on the other side to study the creative effects of the teacher on the resources s/he mobilized. The weaker aspect of the model is to consider the class of situations, and, finally, the situation themselves, as given once for ever. But situations never repeat… Following Bernstein (1996)," a situation is always a repetition without repetition". Vergnaud (2009), in his theory, situates as essential the pair scheme / class of situations, evidencing their joint *evolution* as a key condition of each cognitive process:

> The function of schemes, in the present theory, is both to describe ordinary ways of doing, for situations already mastered, and give hints on how to tackle new situations. Schemes are adaptable resources: they assimilate new situations by accommodating to them. Therefore the definition of schemes must contain ready-made rules, tricks and procedures that have been shaped by already mastered situations; but these components should also offer the possibility to adapt to new situations. On the one hand, a scheme is the invariant organization of activity for a certain class of situations; on the other hand, its analytic definition must contain open concepts and possibilities of inference (Vergnaud 2009, p. 88).

This adaptive aspect of schemes opens then a new question, aiming to rethink the DAD initial model (Figure 4) for taking into account the joint evolution scheme / class of situations.

The third level of reflection is related the *structure* of the documentation work. During the past ten years, we have focused on operational invariants, for inferring teacher's knowledge in action. I think that we should give more importance to the rules of action, of gathering information, and of control, as well as to the inferences, that Vergnaud considers as other components of a scheme. We use to look at a resource system as a structured entity; we should look at a documentation work as a structured entity as well. It opens new questions, as: what are the successive stages of such a work; what are the links between each stage (the events/issues triggering the passage from one stage to another stage); is there a model of documentation work characterizing a teacher's profile, or/and a class of situations, or/and a discipline? What about the model of documentation work of two (or more) teachers working together?

Such a reflection has begun in some works, for example (Trouche et al. 2018a Online first). Studying the preparation of a new lesson by two teachers, a 'DAD driven' analysis points out successive stages of their documentation work (analyzing the curriculum, visiting resources from their own resource systems, comparing-selecting some of them, writing a first version of a lesson plan based on the combination of these resources, trying to integrate this lesson plan in a global progression; evidencing some missing resources and looking for them out of their resource systems, but in their 'resource confident zone'…). This paper (ibid.) evidences also that the structure of a teacher's documentation work resulting from a given analysis, even if it follows a common model (Figure 5), is strongly sensitive to the theoretical lens grounding this analysis; in this case, three different lenses are used, DAD, Cultural and Historical Activity Theory (Engeström 2014) and Anthropological Theory of the Didactic (Chevallard 1999). Such a structural analysis of teachers' documentation seems to be promising and should be continued.

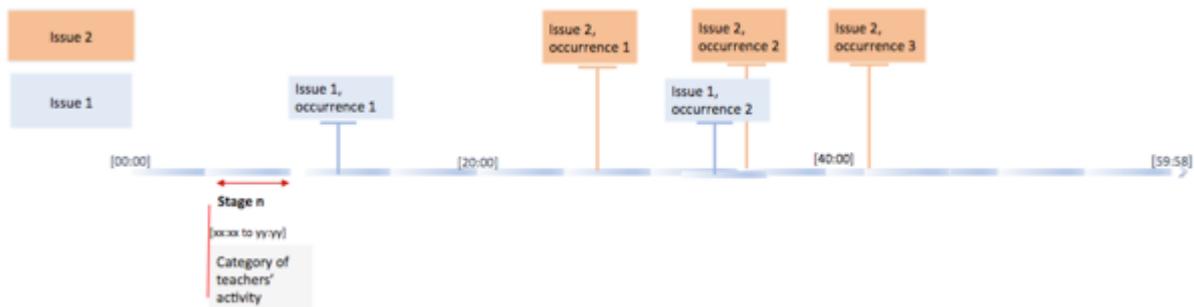

**Figure 5.** *A model for supporting analyses of teachers' documentation work (Trouche et al. 2018a Online first), proposing to cut into successive stages the work of a teacher interacting with resources for a given purpose, and facing different issues*

To me, these three levels (the level of scheme, of situation, and of documentation work) are strongly interrelated, and they motivate my third potential research program: **Deepening the dialectics of schemes / situations of documentation work.**

*3.4 Research program n°4: Deepening the analysis of conditions / effects of teachers collective documentation work*

This program takes into account the second category of students' propositions (§ 3.1), answering the questionnaire, about the dialectic *between individual vs. collective documentation work of teachers*, their context and effects. In this section, we trace the history of this dialectic in the genesis of DAD, before formulating a new research program.

From its beginning (Chapter 1 of this book; Gueudet and Trouche 2008, p. 17), DAD has situated each documentation work as taken in a *bundle of institutional determinations* (from school to society). Thus, among the nine main teachers families of activity that Gueudet and Trouche (2010) distinguished, two are mainly concerned with collectives:

- « Participating in school collective activities (accompanying students for a journey, participating to the school board, monitoring teachers' training);
- Participating in professional collectives out of school (teachers associations, trade unions…). » (p. 68, our translation).

The choice was made, at this starting point, to focus on collectives presenting the strong features of *communities of practice* (Wenger 1998), i.e. a shared commitment, a participation to a shared project, and the existence of a reification process i.e. producing 'things' recognized as a common wealth). The argument was (Gueudet and Trouche 2008, our

Page 13

translation):

> "We retain, in our work, the frame of the communities of practice, because it appears to suit the collectives we want to study, but also because the dialectic participation/reification seems particularly relevant for studying the documentation work. Indeed, it allows to understand the interplay between the commitment in a community and the production of resources" (p. 19).

This choice was indeed productive, allowing us to develop concepts such as *community documentation genesis*, or *community documentation* (Gueudet and Trouche 2008). However, the issue was that most of the collective where teachers meet, sometimes occasionally, and in an informal manner, are not real communities of practice. This was probably the reason why further studies enlarge the scope of the collectives taken into account, calling out other theoretical frames as Hammoud (2012) using the *Cultural and Historical Activity Theory* of Engeström (2004), Rocha (2018) using the theory of *Thought collectives* (Fleck 1934/1981), or Sabra (2011) using the theory of the *Common worlds* (Beguin 2014).

A special issue of ZDM (Pepin et al. 2013a), dedicated to *Re-Sourcing Teacher Work and Interaction*, confirmed the large vision of teachers' collective work, proposing a holistic perspective on collaborative design, (ICT) resources and professional development. It raised new issues of *coherence* and *quality*, coming from teachers sharing resources outside institutional settings. Besides, in this special issue, Gueudet et al. (2013) maintained a focus on communities of practice, retaining three conditions for the development of such communities (p. 1014): a *mutual endeavor*, *minding the system*, and *common forms of addressing and making sense of resources*.

From 2014, new collaboration with China (via the links between ENS de Lyon and East China Normal University) opened a new field of inquiry, giving access to collective work as a regular part of teachers documentation and an essential mean of professional development (Pepin et al. 2016). The Wang's PhD, between France and China, gave means for contrasting the collective documentation work in the two sides (Wang 2018). Chapter 6, in this book, contrasts the collective documentation work between China and Japan.

From 2015 to 2017, a European project, MC2 (see footnote 4), analyzed the creative aspect, for mathematics teachers' documentation work, of combining design processes occurring in community of practices and community of interest (i.e. communities gathering very different people, just sharing the interest for developing creative resources for mathematics teaching). The cross fertilizing of these different communities appeared very productive (Essonier et al. 2018). In the French PREMaTT project (see footnote 2), the objective was to stimulate the design process in a network of schools, considered as 'small factories', supported by researchers and a monthly meeting in a "laboratory for innovative design": then the teachers are exposed to a variety of collectives, in regular as well in artificial settings.

Finally, in the recent period, a focus has been made on "micro-collective." Wang (2018) developed the idea of "documentation-working mates", standing for a pair of teachers working together regularly. This smallest collective appears to be a good frame for deepening the analysis of the collective documentation work, the interactions between the two teachers evidencing some aspects of their schemes (Trouche et al. 2018a Online first).

Then it seems that we have accumulated, in the recent years, a rich experience about teachers collective documentation work. This diversity of collectives, as well as the diversity of theoretical frameworks mobilized, opens new questions, aiming to rethink the DAD initial model (Figure 6) conceived in the case of a Community of Practice. To what extent is it possible to speak of a 'shared repertoire'? Of a collective resource system? To what extent is it possible to speak of common components of scheme? Of shared knowledge? What about teachers belonging to different collectives? How is it possible to take into account the diversity of a collective a given teacher belongs to?

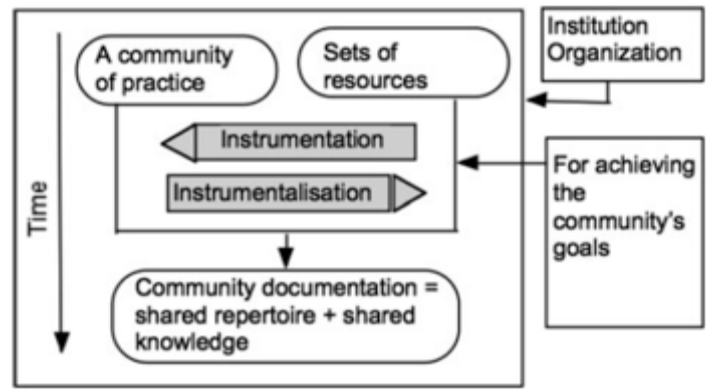

**Figure 6.** *Model for of a collective documentation work in the case of a community* (Gueudet and Trouche 2012b, p. 308)

This set of questions motivates a fourth potential research program: **Deepening the analysis of conditions / effects of teachers collective documentation work, towards an updated model.**

*3.5 Research program n°5: Modeling teachers' working with resources trajectories and professional development over the time*

This program takes into account the fourth category of students' propositions (§ 3.1), answering the questionnaire, about *the length of teachers' professional development* (dialectic between short vs. long geneses). Deepening teachers' documentation work needs to have a large spatial view (all the contexts where this work happens, it is the purpose of the two previous potential programs), and a large time view, it is the purpose of this program. As stated by Pastré (2005, our translation): "There are two main poles of human activity, the first one is structured around the couple scheme – situation […] the second one around the experience acquired, and constantly transformed by the actors" (p. 231). We address in the following discussion: this general issue of time, as an essential component of documentational geneses; the different perspectives allowing us to take into account the short vs. the long periods of time; and, finally, the need for a model of teacher professional development when interacting with resources.

For Vergnaud (2009), time is an essential component of conceptualization: "For instance the analysis of additive structures shows that the concepts of addition and subtraction develop over a long period of time, through situations calling for theorems of very different levels" (p. 89). Adapting this sentence, we could say that teachers' documentation schemes develop over long periods of time, through situations calling for operational invariants of very different levels. The issue is that, constrained by the duration of a Ph.D. (about 4 years) and research projects (most of the time maximum 4 years), we do not have case studies exceeding this length of time. Unfortunately, starting the construction of DAD, we did not planned a follow-up of a collection of case studies, otherwise we would have today some 10 years long case studies. Finally, the only 10 years documentation cases that we have are… ourselves, as researchers interacting with resources.

We have then a number of cases studying the evolution of teachers' documentation work over a short period of time, this period being generally carefully chosen for its hypothetical potential for provoking changes in teachers documentation work, or/and for revealing the strongest invariants allowing the teachers to face the issues at stake: for instances, a period of curricular change (Rocha, Chapter 9); a period of *documentational incident* (Sabra 2016); or a period of *professional metamorphosis* (Assis, Chapter 8, studying the passage from pre-service to in-service teachers). The French ReVEA program (footnote 3) also allowed to analyze, over four years, some deep evolutions (Trouche et al. 2018), due, in mathematics, to the curriculum evolutions (introduction of programming and algorithmics, needing to integrate new resources, and to modify her/his relationship to the ready-made calculation tools), to the new means for storing and sharing information (Dropbox, Google drive…) and to the development of e-textbooks. Indeed the rapid evolution of the teaching environment, and of the schooling form itself (e.g. flipped classrooms) could lead to the idea that 'it is enough to consider short periods of time for capturing major evolutions'. The long history of integration of tools in education (see for example the integration of Interactive White Boards, Karsenti 2016) makes us aware of the necessity to distinguishing between superficial and quick teaching changes on one side, deep and slow ones on the other sides. It is the reason why we need to combine short and long term follow-up of teachers' documentation work.

The ReVEA program aimed to propose, at a national level, and for a set of contrasted disciplines, an institutional observatory for following, over a long period of time, the evolution of teachers' interactions with resources. Unfortunately, this project did not happen. For taking into account the long period of time, we have actually two main perspectives: firstly organizing a long time follow-up and collecting data lively during the whole period (difficult for the reasons evoked above; technically possible in using new tool for data collection, see § 4.3); secondly resting on teachers' 'documentation memory'. Some in progress projects are following this last perspective: Rocha (2018a, Chapter 9), introducing the concepts of *documentational experience* and *documentational trajectory;* Loisy (Chapter 9), introducing the concept of *professional trajectory*; or Santacruz and Sacristan (2018), introducing the concept of *reflective documentational path*. The operationalization of these concepts lies on a double re-construction: a reflective re-construction, by the teacher him/herself from his/her own experience, drawing on the resources s/he used and designed, the events s/he crossed, and the collectives s/he participated to; and a re-construction by the researcher, using the teacher's reflective view, the effective traces of his/her experiences, and the interactions with the teacher, allowing to dig into these traces. This operationalization often results in complex drawing (Figure 7). Zooming in over a long time trajectory to considerer shorter periods could give the impression of a fractal structure, to be investigated further.

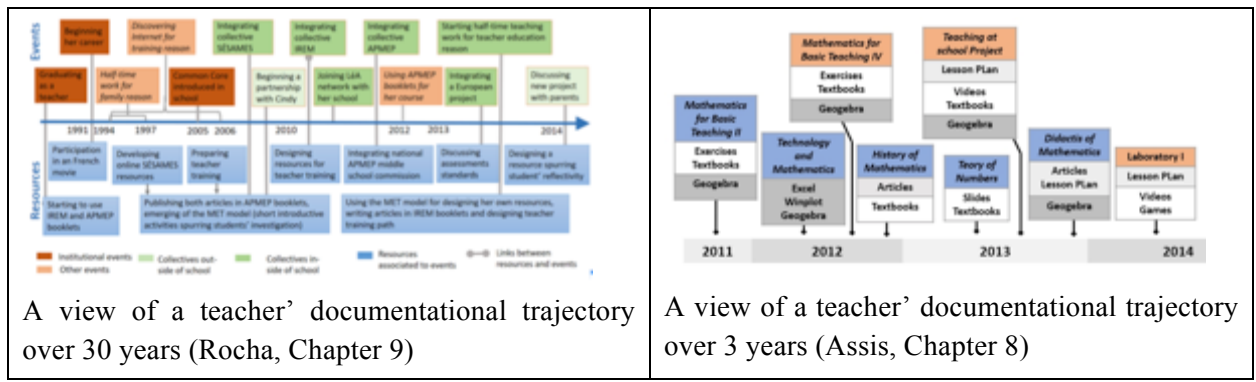

| A view of a teacher' documentational trajectory over 30 years (Rocha, Chapter 9) | A view of a teacher' documentational trajectory over 3 years (Assis, Chapter 8) |

**Figure 7.** *Zooming in, and zooming out of a teacher's documentational trajectory, determined by a sequence of critical events, critical resources used/designed and critical collectives: an appearance of fractal structure*

The analyzes of *teachers' trajectory* (we keep this expression for designating the various concepts introduced above) is necessarily linked to a point of view on teacher professional development in interacting with resources. This professional development has been thought under different models as *Mathematics teaching expertise* (Pepin et al. 2016), *teaching design capacity* (Pepin et al. 2017), *pedagogical design capacity* (Brown 2009), *information competencies* (Messaoui, to be published) or *documentational expertise* (Wang, Chapter 10). This diversity opens a set of questions: how could we collect data allowing us to study teachers' interactions with resources over a long period? Focusing on sensitive periods (as beginning teachers, or changing position)? Which concepts to be introduced for charting the evolution of teachers' experience over this time? How could we combine the study of different trajectories (in various collectives, in various institutions…)? Finally: how can we model teachers' professional development in interaction with resources? (or: To what extent does the diversity of concepts used call up for a multidimensional model of teacher professional development?)

This set of questions motivates my fifth potential research program: **Modeling teachers' working with resources' trajectories and professional development over the time.**

I reviewed in this section, from the missing resources pointed out by the answers to the questionnaire, four potential research programs (numbered 2, 3, 4 and 5). Not surprisingly, concerning a developing theoretical frames, these programs mainly concern concepts and models. I will come back to these programs in a synthetic way later, in the discussion section (§ 5, Figure 15). In the following section, I address some methodological issues.

## 4. Issues of methodological design

I address in this section some methodological issues related to DAD. Starting from the missing resources pointed out by the answers to the questionnaire, I propose then three potential programs of research aiming to take into account these methodological needs.

*4.1 DAD as an incubator of methodology for analyzing teachers' work with resources*

Gueudet (Chapter 1) points out the importance of the development of a specific methodology, called *reflective investigation methodology*, for deepening DAD. I have already evidenced in this Chapter (§ 2.2), from the answers to the questionnaire, to what extent the development of new methodological tools appeared both as a result and a need of/for DAD. Indeed the need for analyzing very heterogeneous resources leads to a very diverse collection of data (see, for instance, Figure 8).

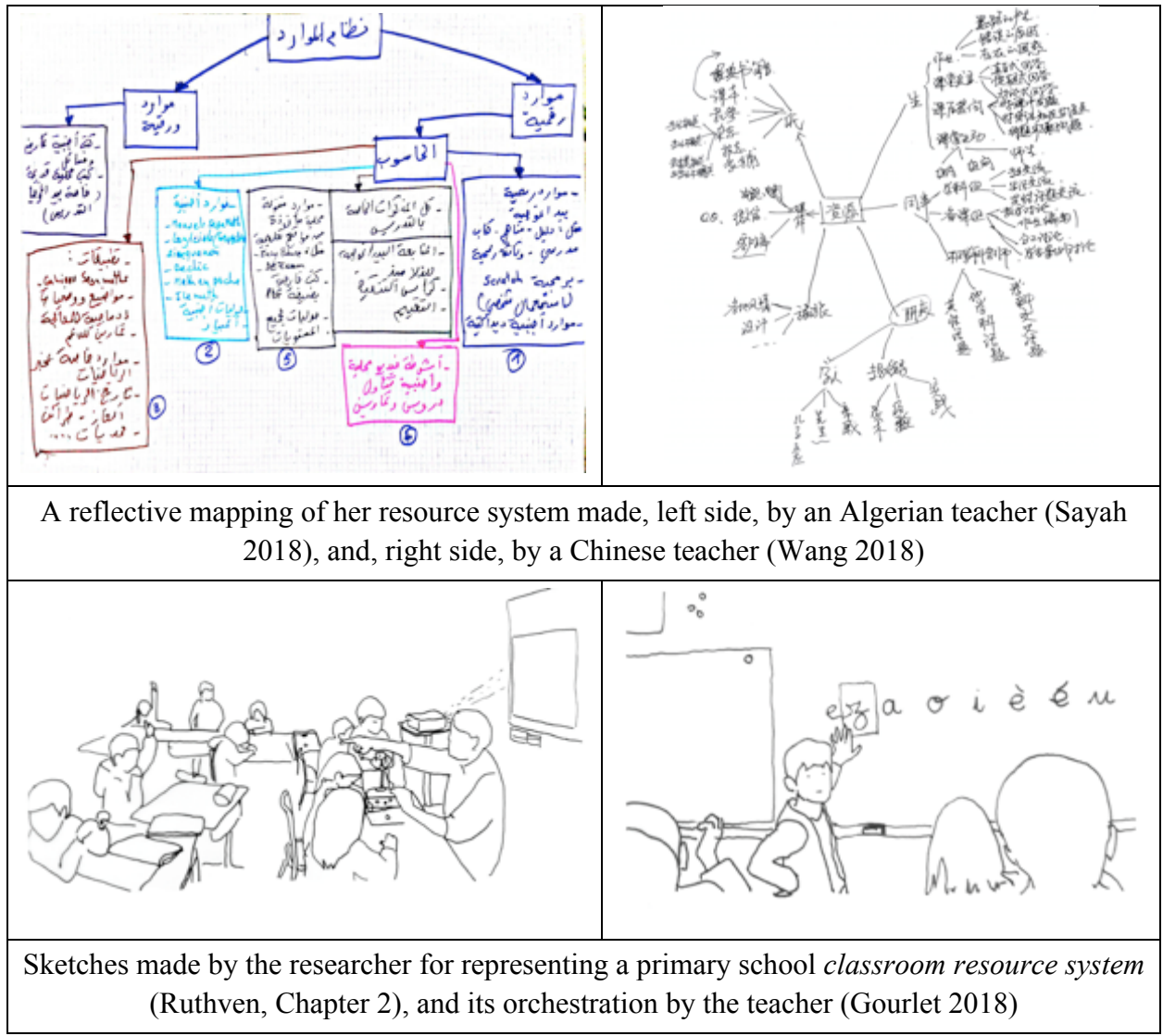

A reflective mapping of her resource system made, left side, by an Algerian teacher (Sayah 2018), and, right side, by a Chinese teacher (Wang 2018)

Sketches made by the researcher for representing a primary school *classroom resource system* (Ruthven, Chapter 2), and its orchestration by the teacher (Gourlet 2018)

**Figure 8.** *Drawings representing resource systems, by the teachers themselves, or by the researcher*

This methodology evolves, as each methodology represents a living theoretical framework. For example, Gueudet, in Chapter 1 (§ 2.2), distinguishes, in her retrospective view of the beginning of DAD, four principles grounding the reflective investigation methodology. Trouche et al. (2018b Online First) add a fifth one: *The principle of permanently confronting the teacher' views on her documentation work, and the materiality of this work (materiality coming, for example, from the collection of material resources; from the teacher's practices in her classrooms)*, evidencing the need to respect the reflective view of the teacher, and the analysis of the researcher.

Chapter 8 in this book, dedicated to methodological issues, evidences both the convergence of the DAD methodological discussions, due to the principles grounding these shared principles, and their diversity, due to the variety of contexts and research questions. This diversity also reveals the underlying complexity of the processes at stake. I have classified answers to the questionnaire into three categories, that could be related to different sensitive points of the methodological model:

- A category on *reflectivity*. How to support and to stimulate it? How to collect and analyze informal data coming from teacher's self productions?
- A category on *tracing resources*. How to trace teachers' interactions? How to study evolving resources? How to manage large scale studies? How to develop quantitative analysis? Which kind of algorithms could be developed for structuring data analysis?
- A category on *representing data*. How could data and data analyses be stored and presented, particularly in the case of collectives? Which kinds of hybrid support could be developed, for sharing the data and discussing the related analyzes, in combining images and sounds, via webdocuments?

Drawing from these categories, I propose in the following three potential research programs.

*4.2 Research program n°6: Looking for methodological models for stimulating reflectivity; storing and analyzing related data*

This program takes into account the first category of students' propositions (§ 4.1), answering the questionnaire, about *stimulating teachers' reflectivity*. Reflectivity is a structural component of the reflective investigation methodology for at least two essential reasons: getting information on teachers' documentation work as a *continuous* process (the need to know what happens between two meetings with the researcher); getting information on documentational geneses as *long* processes (the need to know what happened a long time before). Working on reflectivity raises two issues, not fully addressed until now: how to stimulate reflectivity; how to analyze the very diverse data resulting from the expression of teacher's reflectivity?

*The first issue is to take into account the reflectivity,* that has been done until now mainly through three means: through a *contract* specifying the role of the teacher; through specific *collective settings* for making the teacher confident; through specific situations aiming to stimulate teacher's self investigation. Sabra (2016) has defined what he calls a *methodological contract* as a 'system of *mutual expectations* between the teacher and the researcher': the teacher will describe his/her documentation work, knowing that the researcher is not here for a judgment, or an assessment, but for better understanding 'teachers interacting with resources', producing then new knowledge that could benefit further to the whole profession. Wang (2018) followed two 'documentation-working mates', i.e. teachers preparing regularly their lessons together, and she benefited in this situation of natural reflectivity, each teacher 'speaking to the researcher through her documentation working mate'. Gueudet and Trouche (2010), used a situation proposed by Odonne et al. (1981), the 'Instruction to the double': they asked the teacher to imagine that s/he will be replaced in his/her classroom by another teacher, looking exactly as himself/herself, and that s/he has to transmit his/her own resources and instructions, so that the students would see no difference with their usual teacher and his/her double. These three means (methodological contract, collective settings, instruction to the double) have not been investigated in a systematic way.

I think that we should indeed develop a more systematic 'reflection on reflectivity', its conditions and effects. Vermersch (2012) gives a general frame for rethinking a methodology aiming at the explicitation, by a subject, of his/her action. His book presents the historical, epistemological and practical coherence of the technique of explicitation, as a general issue that had addressed the philosophy over a long time. Drawing from the philosopher Husserl, he evidences the deep links between consciousness, passive memory, and attention. We should draw ourselves from this general frame for rethinking the way we stimulate, and use, teacher's reflectivity about his/her interactions with resources, ancient as well as current.

*The second issue to address is the analysis of data produced by the 'reflective teacher'*, particularly the various representations, or mappings, that a teacher is asked to design, or design for himself/herself. Caraës and Marchand-Zanartu (2011) have evidenced the dynamics of such 'images of thought'. Figure 9 shows the case of Alfred H. Barr mapping his field of interest, reconsidering regularly his view, deepening at each time his knowledge of the field. It gives also an idea of the complexity of the work of transposition – inference, for going from a handmade drawing to a 'proper' printed version. The analysis of such drawings raises a lot of difficulties, which have not been addressed by the diverse cases studies developed under the umbrella of DAD. Hammoud (2012, p. 216) has begun to propose a method, drawing from methods coming from the field of cognitive mapping (Cossette 1994, 2003), which allows to infer some information from teachers' drawings or sketches.

These issues open a wide field of questions, that are, to me, interrelated: Which settings - tools to be thought for stimulating teachers' reflectivity? How to combine, and make profit of, individual, and collective, reflectivity? Which kind of expressions of teacher's reflectivity should we privilege? Which analyses to be develop for understanding the traces of these expressions? And finally, should we look for a single model for stimulating, tracing, analyzing, teachers' reflectivity?

This set of questions motivates my sixth potential research program: **Looking for methodological model(s) for stimulating reflectivity; storing and analyzing related data.**

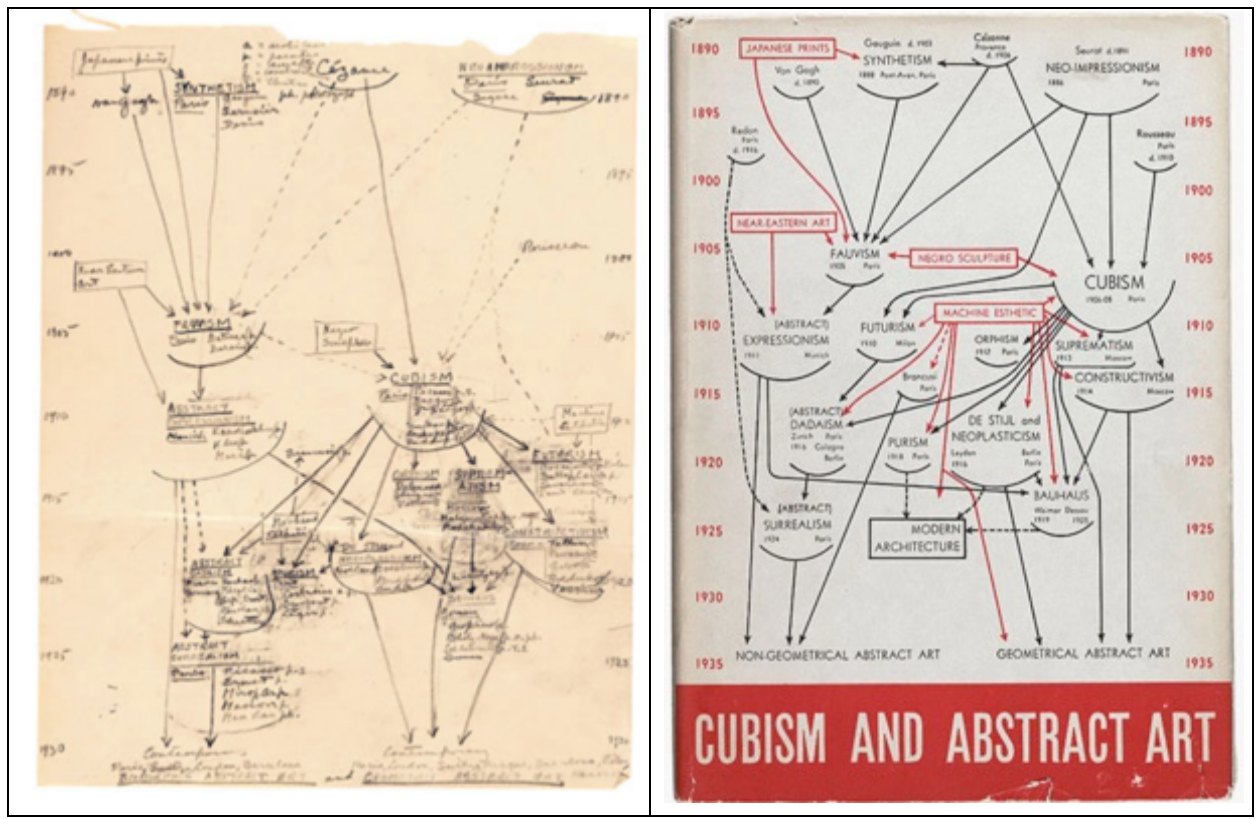

**Figure 9.** Left, one of the successive mappings of Modern Art drawn by Alfred H. Barr, founder of the MOMA *(*Caraës and Marchand-Zanartu 2011, p. 86), and right side, a transposition of these mappings for printing a MOMA poster.

*4.3 Research program n°7: Developing models combining quantitative and qualitative studies of interactions between teachers and resources*

This program takes into account the second category of students' propositions (§ 4.1), answering the questionnaire, about *tracing teachers' resources*. There is, indeed, a contrast between the abundance of data provided by teachers using (individually as well collectively) digital resources, and the frequent 'artisanal' character of methods used by researchers for investigating these data. In this section, I firstly point out some emerging tools for addressing the 'big data' coming from teachers' documentation work; secondly I evidence the remaining complexity, suggesting a new research program.

Sabra analyzed the emails exchanged by a group of 30 teachers, members of the association Sésamath (Sabra 2011, p. 200), designing an online textbook. Just at the time of the discussion about one chapter, he had to take into account about 1000 emails, and this huge number leads him to develop specific methods for analyzing discussion threads. Such quantitative methods remained quite rare in the DAD corpus. At the same time, for analyzing, in the Internet era, the interaction between a large number of users and a large number of resources, new field of research developed, as: *Information Architecture* (Pedauque 2006, already quoted by Gueudet in Chapter 2; Salaün 2012); *Learning Analytics*, for the measurement, collection, analysis and reporting of data about learners and their contexts (already used for education policy purposes, see Ferguson et al. 2016); and Teaching Analytics, as the application of learning analytics techniques to understand teaching and learning processes. Teaching analytics is used, for examples, for analyzing classroom interactions (Prieto et al. 2016) or for analyzing online teachers interactions when they search and create educational resources (Xu and Reker 2012), and developing, in this perspective, specific tools and methods.

Such analytic tools are in progress. Salaün (2016), in a seminar dedicated to the analysis of interactions in a MOOC dedicated to mathematics teachers education, recognized that: "We miss tools for dynamically following up the documentation work of a given community" (Figure 10). Moreover, we need to be aware that it is impossible to totally understand teachers' documentation work by large quantitative studies. Jansen (2018), comparing the efficiency of modeling in experimental and human sciences, underlines that there are four essential factors that make the simulations of society qualitatively more difficult than those of matter: the heterogeneity of humans; the lack of stability of anything; the many relationships to be considered both temporally and spatially; the reflexivity of the humans who react to the models that one makes of their activity.

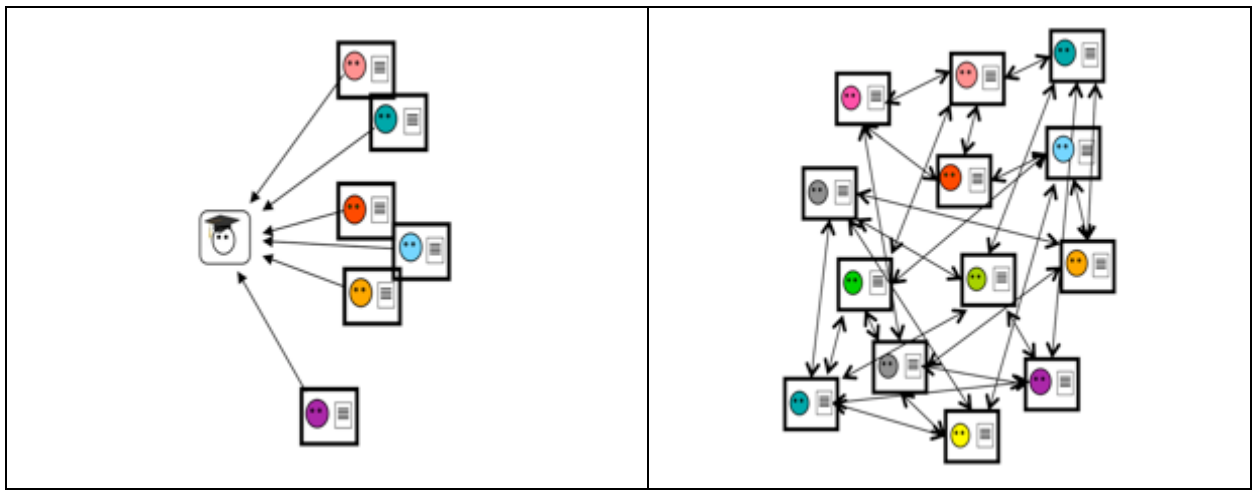

**Figure 10.** *From a teacher centered setting (each student interacting with a teacher) to a collaborative work setting (each learner interacting with each learner) a leap in complexity* (Salaün 2016)

These issues open a set of questions: which interactions could be developed with the field of teaching analytics for modeling teachers' interactions with resources? How to articulate quantitative studies and qualitative case studies? These questions motivate my seventh potential research program: **Developing models combining quantitative and qualitative studies of interactions between teachers and resources.**

*4.4 Research program n°8: Thinking reflective and collaborative supports for capturing, analyzing and sharing data related to teachers' documentation work*

This program takes into account the third category of students' propositions (§ 4.1), answering the questionnaire, about *representing data*. The diversity of interrelated data to be collected when following teachers' documentation work calls up for specific support. In this section, we firstly evoke several attempts aiming to address this issue, then we evidence some main questions needing further research.

A first attempt was the *documentation valise*, described by Gueudet (Chapter 1, § 3.2) A documentation valise gathers a set of data concerning a teacher's case study: the ambition is to make understandable, by a researcher from outside (i.e. not involved in this study), the case and its study. Only one case, Vera's case, is presented on the website dedicated, at its beginning, to DAD[9]. The case focuses on a *lesson cycle*, following four steps (preparing the lesson; implementing it: debriefing throughout assessing students' work; and reflecting the whole process). The interface provides general information on the case's context (curriculum, school, teacher) and on the research context (objectives, methodology). It provides also specific information on the lesson cycle: the videos of each step, the video scripts made by the researcher and the resources used/produced by the teacher. It gives no means to the reader for commenting the data, or proposing an alternative analysis of what happened. Due to this limitation, and as no specific interface was developed for supporting the design of such documentation valises, this methodological perspective was rarely used for studying teachers documentation work, and presented to the community (see Pepin et al. 2015). Finally, the documentation valise remained a metaphor of the documentation work as a journey, needing, as each journey to gather and transport personal resources, and producing new resources (as travelogues…).

After this first attempt, I had the chance to work, in 2015, for two months, with a research team in Brazil, Lematec[10], interested in developing interactive online supports in the context of teacher training, and, from this stay, emerged this idea of *webdocument*, or webdoc (as an abbreviation of 'web documentary') defined by Wikipedia (https://en.wikipedia.org/wiki/Web_documentary) as "a unique medium to create non-linear productions that combine photography, text, audio, video, animation, and infographics based on real time content. This way the publications progresses over several weeks […] the viewer acquires control of navigation, in a way becoming the author or creator of its own personalized documentary". First examples of such webdocuments were developed in the context of this scientific stay (see Lucena and Assis 2015, including English, French and Portuguese versions), allowing to store both data related to a teacher' documentation work, and a preliminary analysis of these data according to a given research question. Compared to the documentation valise, it constitutes a real improvement, constituting a passage from a metaphor to its operationalization. Bellemain and Trouche (2016) describe (French and English versions) the development of this interface, and of the associated reflection.

Following this dynamic, a third attempt occurred in the context of the French ReVEA project, with the AnA.doc platform (Alturkmani et al. to be published). In addition to a glossary (introduced in this Chapter, § 2.3), AnA.doc distinguished two essential levels (Figure 11): a level of situations, for *storing* data, and a level of webdoc, for *analyzing* them. Compared to the Brazilian webdocuments, it could be considered as an improvement, for two points of views: distinguishing these two levels of storage analysis; giving tools for grounding the analysis on excerpts of date (for example one minute video, and not the whole video), leading to more accurate statements. A webdoc could be considered then as a preliminary analysis, a short text with a small number of short excerpts, allowing a quite easy appropriation in order to facilitate the discussion of these preliminary results. Ana.doc is thought as a tool for developing analyses of teachers documentation work within a community of research. Several communications in scientific conferences were then supported by AnA.doc (Messaoui 2018, Rocha 2018b).

---

[9] http://educmath.ens-lyon.fr/Educmath/recherche/approche_documentaire/documentation-valise
[10] http://lematec.net.br/lematecNEW/

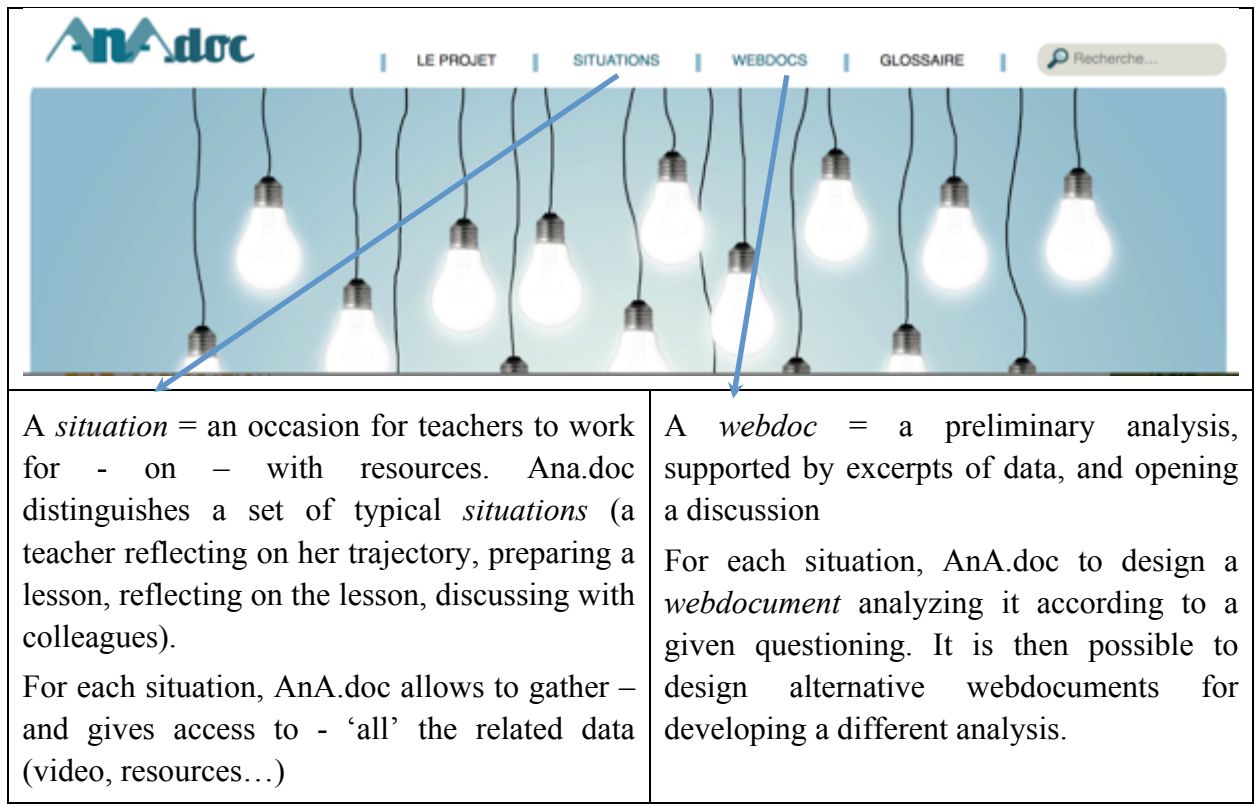

**Figure 11.** *The AnA.doc platform, articulating two essential levels (situations vs. webdocs) for analyzing teachers' documentation work*

A fourth attempt is an ongoing project, PREMaTT (§ 3.4) aiming to develop a twin AnA.doc platform, dedicated to both researchers and teachers, crossing the reflective, collaborative, and research analysis of situations of documentation work.

These developments were studied in the frame of the young researchers workshop of the International Re(s)source conference[11]. Further developments could be then expected. The design of such platform opens actually a set of questions: how could we develop flexible platforms for storing and analyzing data coming from teachers documentation work, that could be used both by small communities and at a large scale? Could we combine the development of such a platform and the combination of individual teachers' travelogues? How could we think travelogues jointly designed and used by teachers and researchers? This set of questions motivates my eighth potential research program: **Thinking reflective and collaborative supports for capturing, analyzing and sharing data related to teachers' documentation work.**

This section was dedicated to the missing resources of DAD, from a methodological point of view, proposing three potential research programs (see § 6, discussion section, Figure 15, for a synthetic presentation). The methodological developments, of course, could not happen from a reflection restricted to DAD, but has to benefit of the interactions with other theoretical frameworks. In the following section, we deepen an issue already addressed in this chapter, the conditions and interest for crossing DAD with other theoretical frameworks.

---

[11] See session C at https://resources-2018.sciencesconf.org/resource/page/id/10

## 5. Towards an extension/expansion of the theory

In this section, I go back to the questionnaire answers regarding the link to be built between DAD and other theoretical frames. Then I examine this question from a general point of view. Drawing from these considerations, I propose finally a ninth potential research program.

*5.1 Back to the questionnaire*

The third point of the questionnaire was (§ 1): *For designing missing resources, which possible theoretical crossings? Which new research programs to be launched?*

In the following, we consider firstly answers questioning the theoretical crossing itself, then we consider the theoretical crossings that are proposed, and, finally, we analyze the current intertwining between DAD and other theoretical frameworks.

There is always, when looking for some theoretical missing resources, a choice to be made between developing the frame housing your research, and looking outside this frame for getting these resources. Some answers to the questionnaire privilege the first choice, proposing, for example, to deepen the concept of mediation in introducing the concept of *documentational mediation* (Sokhna 2018), or to deepen the theory of instrumental orchestration by introducing the concept of *meta-orchestration* (Lucena et al. 2016).

Most of the answers to the questionnaire suggest theoretical crossings: first of all with theories in didactics (didactics of mathematics, joint action in didactics, or professional didactics); secondly with theories speaking of resources and technology (Information and communication sciences, theory of variations, boundary objects, or TPACK); thirdly with theories speaking of collective, history and culture (Lesson studies, Cultural and Historical Activity Theory, Models of Wartofsky, Anthropological approach to didactics, or sociology); fourthly with theories speaking of computation (computer sciences, constructionism and computational thinking).

The need for theoretical means for taking into account the social and cultural aspect is reinforced by answers underlining the necessity to study teachers' work with resources in rural, or poor schools, or specific communities, and the necessity to address the issues of equity and access to resources.

These answers reflect not only the intended crossings with other theoretical frameworks, but existing connections, from the origin of DAD (see Gueudet, Chapter 1; Artigue, Chapter 4; and in this chapter § 3.3 and § 3.4, for the frames mobilized for studying collective aspects of teachers' documentation work). Going further, we could analyze the *resources* produced by the researchers using DAD – their scientific papers – and extract the theoretical references used, for evaluating the degree of connection of DAD with other fields. Tools for making such an analysis exist (Figure 12). We could identify the papers using DAD at an international level, and apply these tools to this set of papers.

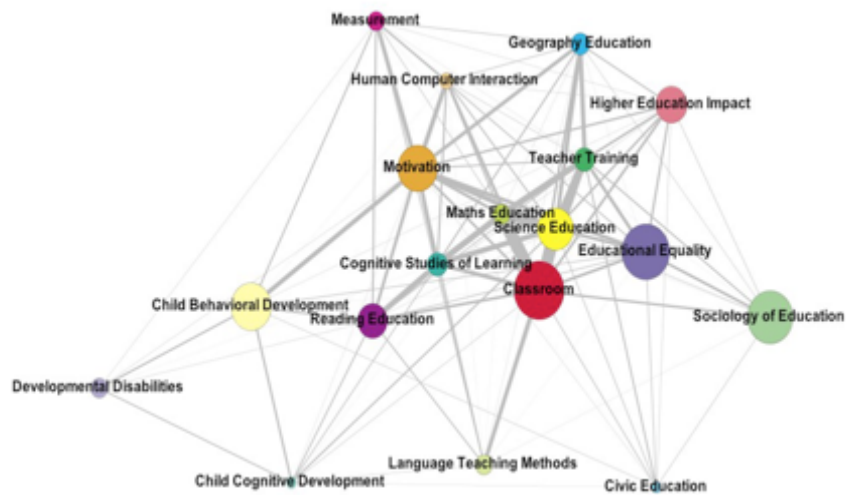

**Figure 12.** *Map of clusters of educational research, from the references used by 22000 papers written between 2000 et 2004 in journals recognized by the French Ministry of Higher Education (see a description of the method in Trouche 2014, p. 2)*

We could develop also this analysis from another point of view: choosing the most 'popular' paper presenting DAD, that is the paper written by Gueudet and Trouche (2009) published by Educational Studies in Mathematics, quoted, according to Google Scholar (on 11[th] January 2019), 287 times; and studying the theoretical frames of the articles quoting this paper.

Beyond these statistical considerations, how could we analyze the interest of connecting DAD with other frames? This is the purpose of the following sub section.

*5.2 Theorizing the interactions of DAD with other scientific frames*

Artigue (Chapter 4) proposes a view of a theoretical frame as a *research praxeology*: "The basic idea is to consider that the model of praxeologies that ATD uses to model human practices might be useful to approach the issue of connection between theories, by making clear that theories emerge from research practices and condition these in return, and that connection between theories involves thus necessarily much more than the theories themselves. They cannot be productively established by just looking for connections between theoretical discourses".

This approach enables her to analyze the development of DAD as a *full research praxeology*, but not a *complete* one, as it continuously extends and diversifies research problématiques, and, in this dynamic, reinforces and diversifies its theoretical connections. It is the case, for example, of a recent paper (Trouche et al. 2018a Online first, already presented § 3.3), using three theoretical lenses for analyzing the work of two teachers preparing together a lesson. This paper is written for a special issue *Curriculum ergonomics* of the International Journal of Educational Research (Choppin et al. 2018 Online first). The issue editors define this new field as "exploring the interaction between the design and use of curriculum materials". Participating to this special issue constitutes indeed on occasion, in facing a new research problématique, to cross DAD with other theoretical frames. To which objectives? To deepen and strengthen DAD as a theoretical construct? To reinforce and diversify its theoretical connections? Or to extend, or expand DAD? I try to address these issues in the following sub section.

*5.3 Research program n°9: Thinking DAD as a theoretical frame interacting with other frames, and as a component of a wider field of research.*

The recent history of DAD is an history of expanding its boarders (see Gueudet, Chapter 2): from mathematics to other disciplines, from secondary schools to pre-primary schools and higher education, and to a variety of contexts (social, cultural, formal / informal), from teachers to a variety of actors of education (including students).

It is perhaps the time, 10 years after its emergence, to question the "raison d'être" of this frame. Taking profit of the acronym DAD as a palindrome, we could translate it either as *documentational approach to didactics*, or as *didactical approach to documentation*. This second 'translation' of the acronym is perhaps more consistent with the research praxeology of DAD, its tasks and techniques, technology and theory. And if DAD defines itself as a didactical approach to documentation, it means of course that there are other approaches to documentation (psychological, sociological…). This is conceptualized in the introduction of the book, proposing the new acronym RAME, standing for "Resources approach to mathematics education", and embracing more largely that DAD this field of research.

Following this thread, for addressing the whole set of disciplines to be taught, I have in mind that we could also introduce the expression *Documentational studies*, as a way of developing *a new multidisciplinary academic field embracing the complexity* (Monteil and Romerio 2017) of "teachers interacting with resources". We could have in mind the example of *Environmental studies* as a "multidisciplinary academic field, which systematically studies human interaction with the environment in the interests of solving complex problems" (Wikipedia).

It opens a set of questions: how could we analyze the position and evolution of DAD with respect to different communities of research, throughout the papers produced by researchers referring to this theoretical frame? How could we define a new research praxeology aiming to study, according to different lenses, the interactions between teachers and resources? This set of questions motivates my ninth potential research program: **Thinking of DAD both as a theoretical frame interacting with other frames, and as a component of a wider field of research.**

From theoretical to social and cultural diversity: the following section tackles this last issue, related to the naming systems used, by teachers, in different contexts.

## 6. Crossing teachers' languages for puzzling out names, meanings and resource systems' structure

In this section, I tackle the issue of taking into account the social and cultural diversity of teachers' documentation work: firstly from the analysis of curriculum, resources and interactions with resources; secondly from the analysis of the naming systems designed-used by teachers over their documentation work; and thirdly for proposing a last potential research program.

*6.1 Contrasting teachers' documentation work across cultural and social boundaries.*

The need for developing comparative studies is mentioned in the answers to the questionnaire. Actually, as soon as we consider teaching resources across the national boundaries, their role in the *'figured world'* of classroom resource systems (Pepin 2009, in the case of England, France and Germany), or as *crucial interfaces between culture, policy and teacher curricular practice* (Pepin et al. 2013b, in the case of France and Norway) appears clearly. Remillard et al. (2014), contrasting the case of Flanders, USA and Sweden, speak of *the voice of curricular material*. Comparative studies seem to be as interesting when the countries are close (see Miyakawa and Xu, Chapter 6 for the cases of China and Japan) as when the countries are far away (Wang 2018 for the cases of China and France). Pepin et al. (Chapter 5) evidence the interest of large comparative studies for deepening the work of teachers as designers. But, till now, the issues of the *words* used by teachers when working with resources has not be really addressed.

*6.2 Teachers' naming systems as revealing the springs of their documentation work*

The issue here is essentially different of that addressed in § 2.3 (giving birth to our first potential research program): it was there the issue of denominating/translating concepts of an emerging scientific field. Here we are speaking of words used by teachers in their daily documentation work, words conceived not as isolated entities, but as *systems* revealing the deep structure of *thought-languages* (Jullien 2015).

I realized the interest of studying these naming systems when teaching in Shanghai a master class, working on interviews of teachers conducting by students themselves (Figure 13). The richness of words showed, for instance, the importance to distinguish the resources designed for oneself, and the resources designed for the school community. Actually, it was impossible to translate from Chinese to English word by word: the only possibility was to analyze the Chinese system of words for inferring the structure of their documentation work, what we did for writing a paper on Chinese teachers' resource system (Pepin et al. 2016).

Category: *The resources designed by the teacher*
- 人格化的材料 the materials processed by the teacher with a strong sense of herself
- 从专业眼光提炼过的材料 materials selected and adapted from a common professional community preference
- 经过教师自己重新理解后的材料 materials re-explained by the teacher
- 转化/活化的材料 transformed / lived materials: the seemingly irrelevant materials developed into academic materials

Category: *Teaching through the way of…*
- 言说 speech, lecture, presentation
- 表演 performance (acting, showing)
- 夸张的肢体动作 exaggerating gestures
- 模拟声响 imitating voices

**Figure 13.** *Analyzing the interview of a teacher (educational philosophy in higher education) by the student who made the interview (personal source, master class on educational research, ECNU, Shanghai 2016)*

This study of teachers naming systems has been already engaged in the frame of the *Lexicon* project (Artigue, *Chapter 4*; Clarke et al. 2017), supported at an international level by a number of researchers. It aims to « document the naming systems (lexicons) employed by different communities speaking different languages to describe the phenomena of the mathematics classroom ».

The phenomena of the mathematics classroom are far from recovering the whole teachers' documentation work. Studying naming systems employed by different communities speaking different languages to describe the phenomena of their documentation work, essentially, remains to be done.

*6.3 Research program n°10: Contrasting naming systems used by teachers in describing their resources and documentation work, towards a deeper analysis of teachers' resource systems*

This study was at the center of the Young Researcher Workshop (session D)[12], following the International Conference Re(s)sources 2018. The session has analyzed data gathered beforehand in different language and cultural contexts (Arabic, Brazilian, Chinese, Dutch, French, Mexican, Turkish and Ukrainian). The data had been collected from interviews, focusing on the preparation of a lesson and on teachers' resource systems. For each of these questions, the teachers have been asked to name their resources, describe their classification and the steps of their documentation work. Preliminary work has been done, describing in English the reality to which these terms refer, and proposing examples and counter examples. The session allowed a cross presentation and analysis of these data (see for example Figure 14 for an Ukrainian case). Understanding a given name needs to study its institutional definition, but also its meaning in a given culture (in Ukraine, the meaning of a 'plan', a 'methodical association'…), for a given teacher, and the place of the corresponding resource in the teacher's documentation work (the position of this particular name in a global naming system). It leads to an anthropological approach to teachers' resources, taking into account linguistic, cultural, social and historical background.

---

**Календарно-тематичний план / Calendar thematic plan**

- A document made by teacher for every grade before every academic year according to the curriculum.
- Identifies the order of themes and lessons to be taught according to the numbers of hours for teaching of the themes specified in the curriculum, numbers of the lessons per week and school schedule.
- Allowed differences with the curriculum up to 10%
- Mandatory document
- Should be examined by the methodical association of teachers of the school, agreed with Vice-director in teaching and upbringing work and approved by Director of the school before the 30th of August.

---

**Figure 14.** *Expliciting a word naming a critical resource for teachers in Ukraine (Rafalska's contribution, see footnote 13)*

This is the beginning of a process: thinking of teachers' resource through the lens of the naming system they develop and use. This process opens a set of questions: Which methodology for gathering relevant data giving access to teachers naming systems? Which interaction with the Lexicon project? How differentiating institutional naming systems, and teachers' personal naming systems? How contrasting naming systems coming from a same culture vs. crossing different cultures (see Wang et al. 2019 for contrasting a Chinese and a Mexican case)? This set of questions motivates my tenth potential research program: **Contrasting naming systems used by teachers in describing their resources and documentation work, towards a deeper analysis of teachers' resource systems.**

---
[12] https://resources-2018.sciencesconf.org/resource/page/id/10

This research program closes the short list of ten potential research programs coming from the answers to my questionnaire… and from my own interpretation and construction.

## 7. Discussion

In this section, I propose firstly a retrospective and critical view on the ten potential research programs; then I review some ongoing works, which could be related to these programs; and finally, I propose a metaphor aiming to capture the complexity of teachers' documentation work.

*7.1 A global view on ten research programs, as ten interconnected perspectives of research and development*

My ambition was not, of course, to determine all the resources that the development of DAD requires, and to infer research programs aiming to design all these resources. I drew these missing resources only from answers I received to the questionnaires, and these answers came from new research problématiques creating new theoretical or methodological needs. The 10 research programs built from these missing resources reflect my own point of view and my 'orchestrating experience'. After the lecture given at the International Re(s)sources conference, when writing this Chapter of this book, I also tried to integrate some inputs of the previous chapters. Even now, I am aware of the incomplete character of this investigation. The search for missing resources looks like a 'mise an abîme', each collection of missing resources appealing for 'finding resources missing in the missing resources that have been found'. Then we could say that it is an imperfect work by nature: impossible to determine a complete missing resource system…

Actually, these ten 'research programs' constitute rather *perspectives* of research and development, as their construction, as programs, remain to be done. Of course, these perspectives are not independent: for example, Program 1 (about a DAD multi-language glossary) is linked with Program 2 (about the structure of teachers' resource systems) and with Program 3 (about documentation schemes). When I look at these ten perspectives of research, I consider finally them as a *kaleidoscope* of perspectives (see Figure 15): according to the point of view that the reader will adopt, when moving around these perspectives, s/he will discover a new panorama of the research to be thought for developing the scientific field.

| | |
|---|---|
| 1. Conceiving a DAD living multi-language glossary | 6. Looking for methodological model(s) for stimulating reflectivity; storing and analyzing related data |
| 2. Modeling the structure and the development of teachers' resource systems | 7. Developing models combining quantitative and qualitative studies of interactions between teachers and resources |
| 3. Deepening the dialectics of schemes / situations of documentation work | 8. Thinking reflective and collaborative supports for capturing, analyzing and sharing data related to teachers' documentation work. |
| 4. Deepening the analysis of conditions / effects of teachers collective documentation work, towards an updated model | 9. Thinking DAD both as a theoretical frame interacting with other frames, and as a component of a wider field of research |
| 5. Modeling 'teachers working with resources' trajectories and professional development over the time. | 10. Contrasting naming systems used by teachers in describing their resources and documentation work, towards a deeper analysis of teachers' resource systems |

**Figure 15.** *A kaleidoscope of 10 perspectives of research*

At this stage of the development of DAD, from these 10 perspectives of research, what stands out is mainly a matter of *names* (names of concepts as well as names used by teachers) and a matter of *models*, what is, in my opinion, the symptoms of a construct coming out of the first experimentations, needing both well defined material and well defined plans for going further. From solid findings (see Gueudet, Chapter 2) to solid foundings...

A last remark on these perspectives. When formulating them, I focused on teachers: *teachers'* resource systems, interactions between *teachers* and resources… Considering the extension of the theory, I think we should consider more widely, beyond the teachers, students, educators… and, finally, each actor involved in a learning/teaching process, these two processes being strongly connected (§ 1.1).

*7.2 Ongoing programs, taking into account some of these perspectives of research*

These perspectives are not proposed in a vacuum. Chapter 1 as well as the current Chapter have shown the variety of works already addressing at least a part of these research perspectives. We could consider the ongoing and future PhDs, the conferences to come, and the research projects.

There are about ten ongoing PhD having chosen the main frame of DAD. I hope that this book, and the research perspectives that it proposes, will inspire other PhDs.

The Re(s)sources 2018 International Conference (Gitirana et al. 2018) had witnessed the vitality of the research in this area. In this book, Chapters 8 to 11 evidence the reflection spurred by the Conference Working groups. The Young Researchers Workshop[13] following the Conference have opened also several research sites. Among them, the reflection and development on didactic metadata for tagging teachers' resources (coordinated by Yerushalmy and Cooper), on webdocuments (coordinated by Bellemain, Alturkmani and Gitirana) and on the naming systems used by teachers in a variety of cultural contexts (coordinated by Rafalska and me) seem to be particularly active.

---

[13] https://resources-2018.sciencesconf.org/resource/page/id/10

The next international conference will take into account these issues, for example the CERME conference (https://cerme11.org/), February 2019, with several Thematic Working Groups mentioning the resources in their titles: Curricular resources and task design in Mathematics Education; Teaching mathematics with technology and other resources; and Teaching mathematics with technology and other resources[14]. Finally, let us quote the ICMT3 (https://tagung.math.uni-paderborn.de/event/1/), September 2019, dedicated to research and developments on textbooks.

Regarding the research to come, I am sure that a number of them will attend the sensitive points that we have underlined, as the issue of collective work: for example, at an international level, the next ICMI study[15] to be launched will be dedicated to Mathematics teachers working and learning in collaborative groups.

Then… A number of improvements to be hoped during the years to come.

*7.3 A fruitful metaphor*

I would like to close this chapter by a metaphor. Editing our first book in French (Gueudet and Trouche 2002), for presenting DAD, we had chosen, on the front page, a Calder's mobile (Figure 16).

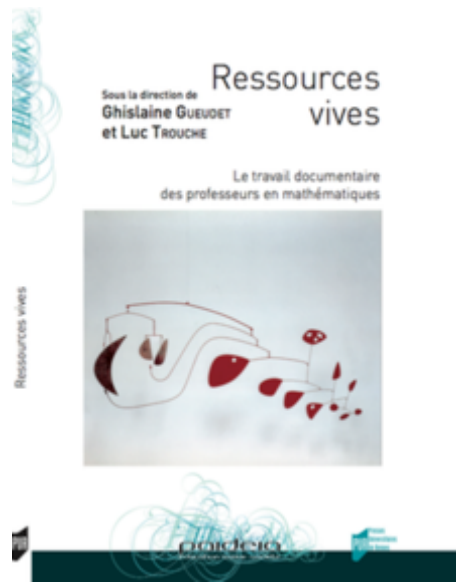

**Figure 16.** *The cover of the first book presenting, in French, DAD (Gueudet and Trouche 2010)*

We had motivated this choice with this short text (p. 371, our translation):

---

[14] However, as we develop in § 6, the reflection on resources is not restricted to the research groups whose title includes the word 'resources'. Taking into account the matter of teachers' work, their *resources*, is certainly a major feature of the current research in educational research.

[15] The ICMI studies are launched by the International Commission on Mathematica Instruction: https://www.mathunion.org/icmi/conferences/icmi-study-conferences. The next study will be the 25th one.

> "A mobile combines dynamics and stability, invariance and movement, like the geneses at the heart of this book; reminds us that "the most profound tendency of all human activity is the march towards equilibrium" (Piaget, 1964). This mobile is part of Calder's series of "gongs", sound mobiles, thanks to collisions between metallic elements produced during their movement; visual and auditory development, visible and invisible resources. It speaks about the individual, interaction, collective ... the necessary commitment of each for a common work".

Some years after, at this stage of development of DAD, and after reading the book of Wohlleben (2017 for the French version) about the secret life of trees, I would choose another metaphor for thinking of teachers' resource systems, as well as to teachers' collective in action: the metaphor of a forest.

> Trees learn, support, interact within a complex, highly interconnected ecosystem, where diversity is an asset, and experience a valued richness. A symbiosis is formed, underground, between tree roots and mycorrhizal fungi. Through this network, carbon, phosphorus, nitrogen or hydrogen, circulate from one tree to another. After cutting a tree, the subterranean relationships between the tree's roots and neighboring trees intensify for a whole period of time, as if to capture the experience of the cut-off tree

Rogerio Ignacio and Cibelle Assis, from Brazil, who were in the French Institute of Education for preparing the Re(s)sources International Conference, proposed to push this metaphor further, in evoking a specific Brazilian tree: the anacardier of Pirangi. Its main feature is to develop a whole forest from one tree, its roots diving in, and emerging from, the earth for developing new trees. This is certainly a fruitful metaphor, as it invited to consider, for understanding the development of a given resource, the whole set of roots feeding it, and linking it to a complex resource system.

Living and supportive resources…

## 8. Acknowledgements